\documentclass[11pt,reqno]{amsart}
\usepackage[T1]{fontenc}
\usepackage{amsmath,amssymb}
\usepackage[margin=1in]{geometry}
\usepackage[hidelinks]{hyperref}

\newtheorem{theorem}{Theorem}[section]
\newtheorem{lemma}[theorem]{Lemma}
\newtheorem{proposition}[theorem]{Proposition}
\newtheorem{corollary}[theorem]{Corollary}
\theoremstyle{definition}
\newtheorem{definition}[theorem]{Definition}
\theoremstyle{remark}
\newtheorem{remark}[theorem]{Remark}
\newtheorem{example}[theorem]{Example}
\DeclareMathOperator{\Ann}{Ann}
\title[A semiprime left ideal]{A semiprime left ideal that is not an intersection\\
of prime left ideals}
\hypersetup{pdftitle={A semiprime left ideal that is not an intersection of prime left ideals},
 pdfauthor={Jakob Cimpri\v{c}}}
\author{Jakob Cimpri\v{c}}
\address{University of Ljubljana, Faculty of Mathematics and Physics,
Department of Mathematics,\newline Jadranska 21, SI-1000 Ljubljana, Slovenia}
\address{Institute of Mathematics, Physics and Mechanics,\newline
Jadranska 19, SI-1000 Ljubljana, Slovenia}
\email{jaka.cimpric@fmf.uni-lj.si}
\thanks{The author acknowledges the financial support from the Slovenian Research and Innovation Agency (ARIS), research core funding No.\ P1-0222 and research projects No.\ J1-60011 and No.\ J1-50002.}
\date{}

\begin{document}

\begin{abstract}
In the free algebra on any uncountable set over any field, we construct
a semiprime left ideal \(I\) and a free generator \(t\) such that \(t\)
lies outside \(I\) but in every prime left ideal containing \(I\).
We also prove that in every countably generated algebra over a
commutative ring, every semiprime left ideal is an intersection of
prime left ideals. Thus a free algebra over a nonzero commutative ring
has this intersection property if and only if its set of generators
is at most countable. Finally, we compare the intersection properties
associated with four other notions of prime left ideal.
\end{abstract}

\maketitle

\section{Introduction}

In any ring, every semiprime two-sided ideal is an intersection of
prime two-sided ideals. We call a ring \emph{prime-separated} if every
semiprime left ideal is an intersection of prime left ideals.
PI rings, rings satisfying the ascending chain condition on semiprime
two-sided ideals, and biregular rings are prime-separated
\cite{Cimpric}. We show that prime separation can fail in a free
associative algebra over a field.

All rings and modules are unital. Modules are left modules, and scalar
rings act centrally. Countable includes finite. A left ideal \(I\) of a ring \(A\) is \emph{semiprime} if
\(aAa\subseteq I\) implies \(a\in I\). A proper left ideal \(P\)
is \emph{prime} if \(aAb\subseteq P\) implies \(a\in P\) or
\(b\in P\). In these terms, \(A\) is prime-separated exactly when
each element \(a\) outside a semiprime left ideal \(I\) is omitted
by some prime left ideal containing \(I\).

Dauns \cite[1.3]{DaunsIdeals} and Hansen \cite[Section~1]{Hansen}
used the right-sided forms of these definitions. Hansen proved the
intersection theorem for right weakly regular rings
\cite[Lemma~2.1]{Hansen}.

\begin{theorem}\label{main}
Let \(k\) be a field and \(X\) an uncountable set. In the free algebra
\(A=k\langle X\rangle\) there exist a semiprime left ideal \(I\) and a
free generator \(t\notin I\) such that every prime left ideal
containing \(I\) contains \(t\).
\end{theorem}

For \(I\subseteq P\), the class \(a+I\) has nonzero image under
\(A/I\to A/P\) exactly when \(a\notin P\). This suggests asking
which nonzero elements of a semiprime module survive in a prime
quotient. For countably generated algebras, every such element does.

An \(A\)-module \(M\) is \emph{semiprime} if
\[
 JM\cap\{m\in M:Jm=0\}=0
 \qquad\text{for every two-sided ideal \(J\) of \(A\).}
\]
Equivalently, \(L\cap\Ann_A(L)M=0\) for every submodule \(L\) of \(M\),
as in \cite[Definition~2.1]{Cimpric}. To deduce this from the ideal
condition, take \(J=\Ann_A(L)\). For the reverse implication, take
\(L=\{m\in M:Jm=0\}\), so \(J\subseteq\Ann_A(L)\).
A nonzero module is \emph{prime} if every element of \(A\) that kills
a nonzero submodule kills the whole module. A submodule is semiprime
or prime if its quotient has that property. For \(M=A\) these
definitions agree with those for left ideals (Lemma~\ref{lem:cyclic}).

\begin{theorem}\label{thm:separation}
Let \(A\) be a countably generated algebra over a commutative ring
\(k\). Every nonzero element of a semiprime \(A\)-module has nonzero
image in some prime quotient.
\end{theorem}

\begin{corollary}\label{cor:separation}
Let \(A\) be as in Theorem~\ref{thm:separation}. Every semiprime
submodule of an \(A\)-module is an intersection of prime submodules.
In particular, every semiprime left ideal of \(A\) is an intersection
of prime left ideals.
\end{corollary}

For a free algebra on finitely many generators over a field,
\cite[Theorem~D]{Cimpric} gives finite direct-sum decompositions into
prime modules and separation by finite-dimensional simple quotients
for finitely presented semiprime modules. Theorem~\ref{thm:separation} treats
arbitrary semiprime modules over countably generated algebras over a
commutative ring.

Over a field, Theorem~\ref{main} and Corollary~\ref{cor:separation}
give the precise countability threshold for free algebras. Passing
to a residue field gives the same threshold over any nonzero
commutative ring (Section~\ref{sec:threshold}).

\begin{corollary}\label{dichotomy}
Let \(k\) be a nonzero commutative ring and \(X\) a set. The free
algebra \(k\langle X\rangle\) is prime-separated if and only if \(X\)
is at most countable.
\end{corollary}

Section~\ref{sec:historical} compares prime left ideals with weakly
prime, extremely prime, completely prime and classically prime left
ideals. Definition~\ref{def:types} gives the corresponding
semiprimeness conditions. Van der Walt proved that a left ideal \(I\)
disjoint from an \(m\)-system is contained in a weakly prime left ideal
disjoint from the same system \cite[Proposition~1.4]{VanDerWalt}.
Here an \(m\)-system is a nonempty subset \(S\) of a ring \(A\) such that
for every \(s,t\in S\) some \(sxt\) with \(x\in A\) belongs to \(S\).

If \(I\) is semiprime and \(a_0\notin I\), one can recursively choose
\(a_{n+1}\in a_nAa_n\setminus I\). To see that \(\{a_n:n\ge0\}\) is
an \(m\)-system, take \(n=\max(i,j)\). Each later term is both a left
and a right multiple of every earlier term, so write \(a_n=a_i u=v a_j\).
Then \(a_{n+1}=a_n r a_n=a_i(urv)a_j\) for some \(r\in A\).
Since every \(a_n\) lies outside \(I\), van der Walt's result gives a
weakly prime left ideal containing \(I\) and omitting \(a_0\).
Proposition~\ref{prop:weak-separation}
proves the corresponding intersection theorem for weakly semiprime
left ideals. Example~\ref{ex:triangular} shows that even two weakly
prime left ideals can have an intersection that fails weak
semiprimeness. Table~\ref{tab:tests} records both implications for
each of the five notions.

\section{The counterexample}\label{sec:counterexample}

In this section \(k\) is a field and \(X\) is uncountable.
The free generators will be a letter \(t\), letters \(b_\alpha\)
indexed by labels, and marker letters \(c_{n,E}\) indexed by a
positive integer \(n\) and a finite set \(E\) of labels.
The defining relations place
\(b_\alpha w b_\alpha\) and \((t-b_\alpha)w(t-b_\alpha)\) in \(I\)
whenever \(w\) has no marker whose finite set contains \(\alpha\).
The additional ideal \(K\) will ensure that each prime left ideal
above \(I\) retains at most one marker for each integer \(n\).

For semiprimeness, a fresh marker whose finite set contains all
labels in a chosen polynomial leaves enough freedom to define a
linear action detecting that polynomial. We give the construction
and both arguments in full.

\subsection{The left ideal}

Take a set \(\Lambda\) of cardinality \(|X|\). Write
\([\Lambda]^{<\omega}\) for the set of its finite subsets, including
the empty set. This set has cardinality \(|X|\). So we can label the
free generators of \(A=k\langle X\rangle\) as follows.
\[
 t,\qquad b_\alpha\quad(\alpha\in\Lambda),\qquad
 c_{n,E}\quad(n\ge1,\ E\in[\Lambda]^{<\omega}).
\]
The letters \(c_{n,E}\) are called \emph{markers}. A word is
\emph{compatible} if whenever it contains \(c_{n,E}\) and \(c_{n,D}\)
for the same \(n\), we have \(E=D\). Let \(K\) be the span of the
incompatible words. It is a two-sided ideal, since a word that
contains an incompatible word is incompatible.

Let \(I\) be the left ideal generated by \(K\) and the elements
\begin{equation}\label{relations}
 b_\alpha w b_\alpha,\qquad (t-b_\alpha)w(t-b_\alpha)
\end{equation}
for all \(\alpha\in\Lambda\) and all words \(w\) containing no marker
\(c_{n,E}\) with \(\alpha\in E\). The empty word is allowed. We say
that these generators of form~\eqref{relations} \emph{belong to} \(\alpha\) and
have \emph{middle word} \(w\).

\begin{lemma}\label{t-outside}
The generator \(t\) does not lie in \(I\), and \(I\) is contained in
a prime left ideal.
\end{lemma}

\begin{proof}
An incompatible word contains two different markers, so its length
is at least two. A generator of form~\eqref{relations} with middle word \(w\)
is a linear combination of words of length \(|w|+2\ge2\). The left
ideal \(I\) is spanned by the incompatible words and by the products
\(vg\) of a word \(v\) and a generator \(g\) of form~\eqref{relations}. Hence every element of \(I\) is a linear
combination of words of length at least two. Since \(t\) is a word of
length one, \(t\notin I\). For the same reason \(I\) lies in the
kernel \(\mathfrak a\) of the homomorphism \(A\to k\) sending every
free generator to zero. The left ideal \(\mathfrak a\) is
prime: if \(a,b\notin\mathfrak a\), then \(ab\notin\mathfrak a\).
\end{proof}

\subsection{Prime left ideals above \texorpdfstring{\(I\)}{I}}

\begin{lemma}\label{core}
Let \(P\) be a prime left ideal of a ring \(A\). Then
\(\mathfrak p=\{u\in A:uA\subseteq P\}\) is the largest two-sided
ideal contained in \(P\). If \(u,v\in A\) and
\(uAv\subseteq\mathfrak p\), then \(u\in\mathfrak p\) or
\(v\in\mathfrak p\).
\end{lemma}

\begin{proof}
The first assertion is clear. Let \(uAv\subseteq\mathfrak p\) and
\(v\notin\mathfrak p\). Take \(s\in A\) with \(vs\notin P\). For each
\(r\in A\) we have \((ur)A(vs)\subseteq uAvA\subseteq P\). Primeness
of \(P\) gives \(ur\in P\) for all \(r\), hence \(u\in\mathfrak p\).
\end{proof}

\begin{proposition}\label{t-inside}
Every prime left ideal \(P\) of \(A=k\langle X\rangle\) that contains
\(I\) contains \(t\).
\end{proposition}

\begin{proof}
Let \(\mathfrak p\) be as in Lemma~\ref{core}. Since \(K\) is a
two-sided ideal contained in \(P\), we have \(K\subseteq\mathfrak p\).

Fix \(n\). If \(E\ne D\), every word in \(c_{n,E}Ac_{n,D}\) is
incompatible, so \(c_{n,E}Ac_{n,D}\subseteq K\subseteq\mathfrak p\),
and Lemma~\ref{core} gives \(c_{n,E}\in\mathfrak p\) or
\(c_{n,D}\in\mathfrak p\). Hence at most one marker \(c_{n,E}\) lies
outside \(\mathfrak p\).

Let \(\Lambda_0\) be the union of the finite sets \(E\) for which
some marker \(c_{n,E}\) lies outside \(\mathfrak p\). There is at most
one such set for each \(n\), so \(\Lambda_0\) is countable. As
\(\Lambda\) is uncountable, we can choose
\(\alpha\in\Lambda\setminus\Lambda_0\).

Let \(w\) be a word. Suppose it contains a marker \(c_{n,E}\) with
\(\alpha\in E\). If \(c_{n,E}\notin\mathfrak p\), then
\(E\subseteq\Lambda_0\) by definition of \(\Lambda_0\), contrary to
\(\alpha\notin\Lambda_0\). Thus \(c_{n,E}\in\mathfrak p\). Writing
\(w=u c_{n,E}v\), we get \(w\in\mathfrak p\) since \(\mathfrak p\) is a
two-sided ideal. Then both elements
\(b_\alpha wb_\alpha\) and \((t-b_\alpha)w(t-b_\alpha)\) lie in
\(\mathfrak p\subseteq P\). For every other word \(w\), these two
elements are generators of form~\eqref{relations}. Thus
\[
 b_\alpha Ab_\alpha\subseteq P\qquad\text{and}\qquad
 (t-b_\alpha)A(t-b_\alpha)\subseteq P .
\]
Primeness gives \(b_\alpha\in P\) and \(t-b_\alpha\in P\), hence
\(t\in P\).
\end{proof}

\subsection{Semiprimeness}

We must show that for every \(a\notin I\) there is \(r\in A\) with
\(ara\notin I\). We will take for \(r\) a single marker.

For a word \(u\), let \(\tau(u)\) be the word obtained by deleting all
letters except the markers. We call \(\tau(u)\) the \emph{marker word}
of \(u\) and use it as the degree of \(u\). For each marker word
\(\sigma\), including the empty word, let \(A_\sigma\) be the
\(k\)-span of the words of degree \(\sigma\). The word basis gives
the direct sum below, and \(\tau(uv)=\tau(u)\tau(v)\) gives the
multiplication rule.
\[
 A=\bigoplus_\sigma A_\sigma,\qquad
 A_\sigma A_\varrho\subseteq A_{\sigma\varrho}.
\]
Thus \(A\) is graded by the free monoid of marker words, with
concatenation as the product of degrees.

\begin{lemma}\label{homogeneous}
The left ideal \(I\) is homogeneous for the grading by marker words:
an element of \(A\) lies in \(I\) if and only if each of its
homogeneous components does.
\end{lemma}

\begin{proof}
Each incompatible word is homogeneous. All words occurring in a
generator of form~\eqref{relations} have marker word \(\tau(w)\),
since \(t\) and \(b_\alpha\) are not markers. Left multiplication by
a word \(v\) changes a degree \(\sigma\) to \(\tau(v)\sigma\).
Thus \(I\) is spanned by homogeneous elements. Grouping a linear
combination of these elements by degree shows that each homogeneous
component still belongs to \(I\).
\end{proof}

Next we describe \(I\) on a small alphabet.

\begin{lemma}\label{restriction}
Let \(F\) be a finite set of free generators that contains \(t\) and
contains at most one marker \(c_{n,E}\) for each \(n\). Put
\(A_F=k\langle F\rangle\). Then the generators of form~\eqref{relations} that belong to
a label \(\alpha\) with \(b_\alpha\in F\) and have a middle word in the
alphabet \(F\) generate \(I\cap A_F\) as a left ideal of \(A_F\).
\end{lemma}

\begin{proof}
Let \(J_F\) be the left ideal of \(A_F\) generated by the indicated subset
of the generators of form~\eqref{relations}. Clearly \(J_F\subseteq I\cap A_F\).

Let \(\pi\colon A\to A_F\) be the algebra homomorphism fixing the
letters in \(F\) and sending all other letters to zero. On the
\(k\)-span of the letters, let \(s\) agree with \(\pi\), except that
\(s(b_\alpha)=t\) whenever \(b_\alpha\notin F\). Define a \(k\)-linear
map \(\rho\colon A\to A_F\) on words by
\[
 \rho(1)=1,\qquad \rho(vx)=\pi(v)s(x)
 \quad\text{for a word \(v\) and a letter \(x\).}
\]
Thus \(\rho\) uses \(\pi\) on the prefix and \(s\) on the last
letter. It fixes \(A_F\). It also kills every incompatible word,
because such a word contains a marker outside \(F\), and both
substitutions send that marker to zero. Hence \(\rho(K)=0\).

We check that \(\rho\) sends every left word multiple of a generator of
form~\eqref{relations} into \(J_F\). Write such a generator as \(hwh\),
where \(h=b_\alpha\) or \(h=t-b_\alpha\). Since \(h\) is a
linear combination of letters, the rule for the last letter extends
by linearity to give
\[
 \rho(vhwh)=\pi(vhw)s(h)=\pi(v)\pi(h)\pi(w)s(h)
\]
for every word \(v\). If \(b_\alpha\notin F\), then
\(\pi(b_\alpha)=0\) and \(s(t-b_\alpha)=t-t=0\). Thus the image of
\(vb_\alpha wb_\alpha\) vanishes because of its prefix, and the image
of \(v(t-b_\alpha)w(t-b_\alpha)\) vanishes because of its last factor.
If \(b_\alpha\in F\), both substitutions fix \(h\). A letter outside
\(F\) in \(v\) or \(w\) makes \(\pi(v)\) or \(\pi(w)\) zero.
Otherwise \(\rho(vhwh)=vhwh\), and \(hwh\) is a generator of \(J_F\),
so its left multiple belongs to \(J_F\).

Since \(I\) is spanned by the incompatible words and the left word
multiples of its generators of form~\eqref{relations}, we have \(\rho(I)\subseteq J_F\).
For \(a\in I\cap A_F\), it follows that \(a=\rho(a)\in J_F\).
Thus \(I\cap A_F=J_F\).
\end{proof}

\begin{proposition}\label{semiprime}
The left ideal \(I\) is semiprime.
\end{proposition}

\begin{proof}
Take \(a\notin I\). By Lemma~\ref{homogeneous}, some component \(f\)
of \(a\), with marker word \(\sigma\), lies outside \(I\). The word
\(\sigma\) is compatible, since otherwise every word in \(f\) would
be incompatible and \(f\in K\subseteq I\). Let \(F\) consist of
\(t\) and the letters occurring in \(f\). Choose \(c=c_{n,E}\),
where no marker with index \(n\) occurs in \(a\) and
\(E=\{\alpha:b_\alpha\in F\}\) is finite. Then \(c\notin F\), and
\(F\cup\{c\}\) contains at most one marker for each index.

Put \(B=k\langle F\rangle\), \(J=I\cap B\), \(V=B/J\) and \(\xi=1+J\).
Then \(f\xi\ne0\). Choose a linear operator on \(V\) that sends
\(f\xi\) to \(\xi\), and let \(c\) act by this operator. This makes
\(V\) a module over the free algebra \(B'=k\langle F,c\rangle\).

We claim that \(I\cap B'\) kills \(\xi\). By Lemma~\ref{restriction}
applied to \(F\cup\{c\}\), it suffices to check the generators of
form~\eqref{relations} that belong to a label \(\alpha\) with \(b_\alpha\in F\)
and have a middle word in the alphabet \(F\cup\{c\}\). Every such
\(\alpha\) lies in \(E\), so the middle word cannot contain \(c\). Thus each
such generator already lies in \(J\) and kills \(\xi\), as do
all its left multiples.

Now \(fcf\xi=f(c(f\xi))=f\xi\ne0\), so \(fcf\notin I\). For words
\(u,v\) occurring in \(a\), the marker word of \(ucv\) is
\(\tau(u)c\tau(v)\). Since \(c\) occurs in neither \(u\) nor \(v\),
this equals \(\sigma c\sigma\) exactly when
\(\tau(u)=\tau(v)=\sigma\). Now \(f\) is precisely the linear combination
of the terms of \(a\) with marker word \(\sigma\).
Thus the component of \(aca\) with marker word
\(\sigma c\sigma\) is \(fcf\). If \(aca\in I\), then
Lemma~\ref{homogeneous} would give \(fcf\in I\), a contradiction.
Hence \(aca\notin I\).
\end{proof}

Theorem~\ref{main} follows from Lemma~\ref{t-outside} and
Propositions~\ref{t-inside} and~\ref{semiprime}.

\subsection{Remarks on the construction}

\begin{remark}[Where uncountability is used]
It is used once, in the proof of Proposition~\ref{t-inside}, to
choose a label outside the countable set \(\Lambda_0\). The
construction of \(I\) and the proof of semiprimeness work for every
set \(\Lambda\). For countable \(\Lambda\) the marker sets
kept by a prime left ideal can cover \(\Lambda\), in agreement with
Theorem~\ref{thm:separation}.
\end{remark}

\begin{remark}[The two indices of a marker]
The finite set \(E\) is needed for semiprimeness: one marker must
protect all letters \(b_\alpha\) of a given polynomial at once. The
integer \(n\) provides, for every element \(a\), a marker that does
not occur in \(a\) and is compatible with the markers of \(a\). The
ideal \(K\) limits a prime left ideal to one marker for each \(n\).
\end{remark}

\begin{remark}[Left and right]
The construction is symmetric under reversal of words, except that
\(I\) is a left ideal. The right ideal generated by the same elements
is a semiprime right ideal with the corresponding property.
\end{remark}

\section{Prime separation for countably generated algebras}\label{sec:countable}

We first check the relation between the module definitions and the
left-ideal definitions used in Section~\ref{sec:counterexample}.

\begin{lemma}\label{lem:cyclic}
For a left ideal \(L\) of a ring \(A\), the module \(A/L\) is
semiprime if and only if \(L\) is semiprime. It is prime if and only
if \(L\) is prime.
\end{lemma}

\begin{proof}
Suppose \(A/L\) is semiprime and \(aAa\subseteq L\).
For \(J=AaA\), the class \(a+L\) lies in \(J(A/L)\) and is killed
by \(J\). Hence \(a\in L\). Conversely, suppose \(L\) is semiprime.
For a two-sided ideal \(J\), we have \(J(A/L)=(J+L)/L\).
A class in this submodule killed by \(J\) has the form \(j+L\),
with \(j\in J\) and \(Jj\subseteq L\). Then \(jAj\subseteq L\), so
\(j\in L\). Thus \(A/L\) is semiprime.

Suppose \(A/L\) is prime and \(aAb\subseteq L\), with \(b\notin L\).
The element \(a\) kills the nonzero cyclic submodule \(A(b+L)\), so
it kills \(A/L\), and \(a\in L\). Conversely, suppose \(L\) is prime
and \(a\) kills a nonzero submodule \(H\) of \(A/L\).
Choose \(b+L\ne0\) in \(H\). For every \(r\in A\), we have
\((ar)Ab\subseteq L\), so \(ar\in L\). Thus \(a\) kills \(A/L\),
as required.
\end{proof}

We now prove Theorem~\ref{thm:separation}. Here \(k\) is an arbitrary
commutative ring.

We may work with \(A=k\langle X\rangle\), the free algebra on a
countable set \(X\). Indeed, there is a surjection from this free
algebra onto the given algebra. For modules over the given algebra,
restriction along this surjection leaves the submodules unchanged. Two-sided ideals act through their images, and
every two-sided ideal of the quotient is the image of its preimage.
Consequently semiprimeness and primeness agree over the two rings.

Fix a semiprime \(A\)-module \(M\) and \(0\ne q\in M\).
We first obtain a prime quotient over a countable subring.
Then we imitate finite pieces of the original action to map \(M\)
into an ultrapower. We test primeness there on the cyclic submodule
generated by the image of \(q\). A failure of that test would already
occur at a single coordinate, in the countable prime quotient.

\subsection{A prime quotient over a countable ring}

\begin{lemma}\label{lem:countable}
Let \(B\) be a countable ring and \(U\) a semiprime \(B\)-module.
For each \(0\ne v\in U\), there is a prime quotient of \(U\) in which
\(v\) survives.
\end{lemma}

\begin{proof}
List the elements of \(B\) as \(b_0,b_1,\ldots\), each infinitely often.
Start at \(v_0=v\). If \(b_nBv_n\ne0\), choose
\(0\ne v_{n+1}\in b_nBv_n\). Otherwise record \(b_n\) and keep
\(v_{n+1}=v_n\). Every \(v_n\) is nonzero and the cyclic submodules
\(Bv_n\) form a decreasing sequence. A recorded \(b_n\) kills
\(Bv_n\), so the ideal \(Bb_nB\) also kills \(Bv_n\).

Let \(J\) be the ideal generated by the recorded elements. No \(v_n\)
lies in \(JU\). Indeed, membership would already hold in \(IU\)
for an ideal \(I\) generated by finitely many recorded elements.
For sufficiently large \(m\), all those elements have been recorded
and \(m\ge n\). Then \(Iv_m=0\) and \(v_m\in Bv_n\subseteq IU\),
contrary to semiprimeness.

By Zorn's lemma choose a submodule \(N\) maximal among those containing
\(JU\) and avoiding every \(v_n\). A chain of such submodules has a
union with the same properties. Since \(JU\subseteq N\),
\(\eta(z+N)=\eta z+N=0\) for all \(\eta\in J\) and \(z\in U\).
Thus \(J\) annihilates \(U/N\), so \(J\subseteq\Ann_B(L)\) for every
submodule \(L\subseteq U/N\).

Let \(L\) be a nonzero submodule of \(U/N\). Its preimage properly
contains \(N\), so maximality gives \(v_n+N\in L\) for some \(n\).
Suppose \(b\) kills \(L\) and \(b\notin J\). Choose \(m\ge n\) with
\(b_m=b\). Since \(J\) contains every recorded element, the rule at
this occurrence gives \(v_{m+1}=bcv_m\) for some \(c\in B\).
As \(v_m\in Bv_n\), we have \(cv_m+N\in L\). Applying \(b\) gives
\(v_{m+1}\in N\), a contradiction. Hence \(\Ann_B(L)=J\) for every
nonzero \(L\), so \(U/N\) is prime and \(v+N\ne0\).
\end{proof}

\subsection{Finite tables and the countable subring}

A \emph{selection} consists of finitely many scalars, vectors and
letters, including the scalar \(1\) and the vector \(q\). Its
\emph{table} lists all true equations of the following forms among
the selected entries, with \(s,t,h\in k\), \(y,z,w\in M\) and \(x\in X\).
\[
 s+t=h,\qquad st=h,\qquad y+z=w,\qquad sy=z,\qquad xy=z .
\]
A \emph{reproduction} in a unital subring \(k'\subseteq k\) assigns a
value in \(k'\) to each selected scalar and a value in \(M\) to each
selected vector, keeping \(1\), \(q\) and the letters fixed and
satisfying all the equations in the table.

\begin{lemma}\label{lem:tables}
There is a countable unital subring \(k_0\subseteq k\) in which every
table has a reproduction. For \(A_0=k_0\langle X\rangle\), there is a
prime \(A_0\)-quotient \(\pi\colon M\to V\) with \(p=\pi(q)\ne0\).
\end{lemma}

\begin{proof}
Number the selected scalars and vectors separately and replace each
with the corresponding unknown. Keep the letters fixed and prescribe
\(1\) and \(q\) for their respective unknowns. There are countably many
numbered systems, since \(X\) is countable. Choose one solution in
\(k\) and \(M\) for each system arising from a table, and let \(k_0\)
be the unital subring generated by all the scalar values in these
solutions. It is countable, and the chosen solutions reproduce every
table in \(k_0\).

Put \(A_0=k_0\langle X\rangle\). Every word in \(X\) belongs to \(A_0\),
so \(A=kA_0\). The ring \(A_0\) is countable, and \(M\) is semiprime
over \(A_0\). Indeed, if \(I\) is a two-sided ideal of \(A_0\),
then \(kI\) is a two-sided ideal of \(A\) because \(k\) is central. The ideals \(I\) and \(kI\) have the same
annihilated vectors and satisfy \(IM=(kI)M\), since
\((\lambda a)m=a(\lambda m)\) for \(\lambda\in k\).
Apply Lemma~\ref{lem:countable} to obtain a prime \(A_0\)-quotient
\(\pi\colon M\to V\) with \(p=\pi(q)\ne0\).
\end{proof}

Fix \(k_0,A_0,\pi,V,p\) from Lemma~\ref{lem:tables} for the rest of
the proof.

\subsection{Finitely many words suffice}

\begin{lemma}\label{lem:finite-tests}
Fix words \(w_1,\ldots,w_d\) in \(X\), allowing the empty word.
There is a finite set \(F\) of words with the following property.
If \(a,r\) are in the \(k_0\)-span of \(w_1,\ldots,w_d\),
\(rp\ne0\), and \(a u r p=0\) for every \(u\in F\), then \(aV=0\).
\end{lemma}

\begin{proof}
In \(V\), a scalar from \(k_0\) killing a nonzero vector kills the
cyclic \(A_0\)-submodule it generates. Hence it kills \(V\) by
primeness. It follows that \(V\) is torsion-free over
\(k_0/\Ann_{k_0}(V)\), which is a domain. Indeed, if a product of
two scalars kills \(V\) and the first acts injectively, the second
kills \(V\). Thus \(V\) embeds in the vector space \(V_K\) obtained by
extending scalars to the fraction field \(K\) of that domain.

Identify \(V\) with its image in \(V_K\). For each word \(u\), define
\[
 L_u\colon K^{d^2}\longrightarrow V_K,\qquad
 (t_{j\ell})\longmapsto\sum_{j,\ell} t_{j\ell} w_j u w_\ell p .
\]
Choose a finite set \(F\) of words for which \(\bigcap_{u\in F}\ker L_u\)
has smallest dimension. Intersecting it with any further kernel
leaves it unchanged. Thus it equals the intersection over all words.
The choice of \(F\) depends only on \(V\), \(p\) and the word list.

Now write \(a=\sum_j c_jw_j\) and \(r=\sum_\ell \beta_\ell w_\ell\),
with coefficients in \(k_0\). Read those coefficients in \(K\).
The product \(a u r p\) is \(L_u((c_j\beta_\ell)_{j,\ell})\).
If these products vanish for \(u\in F\), they therefore vanish in
\(V_K\) for every word \(u\), hence in \(V\) by injectivity.
Words span \(A_0\) over \(k_0\), so \(a\) kills \(A_0rp\).
Since \(rp\ne0\) and \(V\) is prime, \(aV=0\).
\end{proof}

\subsection{The ultrapower and the original action}

Index the coordinates by the finite selections \(i\). At each
coordinate fix a reproduction in \(k_0\) by Lemma~\ref{lem:tables}.
Write \(\lambda_i\in k_0\) and \(m_i\in M\) for the replacements of
\(\lambda\in k\) and \(m\in M\), using zero when an entry is absent.
Thus \(1_i=1\) and \(q_i=q\). The scalar \(0\) is always replaced by
\(0\), since \(0+0=0\) belongs to any table that selects it. For
\(a=\sum_j c_jw_j\in A\) put \(a(i)=\sum_j(c_j)_i w_j\in A_0\).

For each fixed selection, take the set of coordinates containing it.
A finite intersection of such sets is the set of coordinates
containing the union of the corresponding selections, so is nonempty.
The supersets of these coordinate sets form a proper filter. This
family contains the full index set, is closed under finite
intersections and supersets, and excludes the empty set. By Zorn's
lemma it extends to a maximal proper filter, an ultrafilter. Adjoining
a set absent from a maximal proper filter would destroy properness,
so the set is disjoint from some member of the filter and its
complement is in the filter.

Call a set \emph{large} if it belongs to this ultrafilter. The
ultrapower \(V^*\) consists of classes \([v_i]\) of families in \(V\),
where families agreeing on a large set are identified. Addition is
coordinatewise. In particular, \([v_i]\ne0\) exactly when \(v_i\ne0\)
on a large set, by the complement property.

\begin{lemma}\label{lem:ultrapower}
The formula \(a[v_i]=[a(i)v_i]\) makes \(V^*\) an \(A\)-module.
The map \(f\colon M\to V^*\) given by \(f(m)=[\pi(m_i)]\) is
\(A\)-linear and sends \(q\) to the nonzero constant class \(p^*=[p]\).
\end{lemma}

\begin{proof}
Define \(\lambda[v_i]=[\lambda_i v_i]\) and \(x[v_i]=[xv_i]\) for
\(\lambda\in k\) and \(x\in X\). These operations respect agreement on a
large set. Replacements can differ from the original scalars even
for \(\lambda\in k_0\), so this action can differ from the coordinatewise
action of \(k_0\). For instance, the equation
\((\lambda\mu)_i=\lambda_i\mu_i\) holds at every coordinate containing
\(\lambda,\mu,\lambda\mu\). Scalar addition works in the same way.
Together with \(1_i=1\), this gives a \(k\)-module structure. Each
letter acts \(k\)-linearly because \(k_0\) is central in \(A_0\).
This defines the free algebra action \(a[v_i]=[a(i)v_i]\).

The table equations for vectors give \(A\)-linearity of \(f\).
For example, on the coordinates containing \(\lambda,m,\lambda m\),
we have \(\pi((\lambda m)_i)=\lambda_i\pi(m_i)\), which gives
\(f(\lambda m)=\lambda f(m)\). Vector addition and the action of each
letter work in the same way. Finally \(q_i=q\) at every coordinate,
so \(f(q)=[p]\ne0\).
\end{proof}

\begin{remark}
If \(k=\mathbb R\), a table containing \(1\), \(2\)
and \(\sqrt{2}\) forces \(k_0\) to contain a square root of \(2\).
A transcendental scalar \(\theta\) outside \(k_0\) nevertheless acts
on \(V^*\) by \([v_i]\mapsto[\theta_i v_i]\), with \(\theta_i\in k_0\)
at each coordinate.
\end{remark}

\subsection{Why the map gives a prime quotient}

\begin{lemma}\label{lem:prime-ultrapower-quotient}
There is an \(A\)-submodule \(G\subseteq V^*\) such that \(V^*/G\)
is prime and \(p^*\notin G\).
\end{lemma}

\begin{proof}
Suppose that \(a\in A\) kills the submodule generated by a nonzero
\(rp^*\), with \(r\in A\). If \(a\) acts nontrivially on \(V^*\),
choose \([y_i]\) with \(a[y_i]\ne0\). Fix a finite word list containing
the supports of \(a\) and \(r\), and take the finite tests \(F\)
from Lemma~\ref{lem:finite-tests}. The same tests work at all coordinates
because the words stay fixed as the coefficients vary.
Each of the following conditions holds on a large set of coordinates:
\[
 r(i)p\ne0,\qquad a(i)y_i\ne0,\qquad a(i)u r(i)p=0
 \quad(u\in F).
\]
The first two follow because the indicated classes are nonzero.
For a word \(u\), Lemma~\ref{lem:ultrapower} gives
\(a u r p^*=[a(i)u r(i)p]\), since \(u\) acts coordinatewise.
These classes vanish for \(u\in F\) because \(a\) kills \(Arp^*\).
There are only finitely many conditions, so they all hold at some
one coordinate. At this coordinate the finite tests give \(a(i)V=0\)
from \(r(i)p\ne0\), contradicting \(a(i)y_i\ne0\).
Thus any \(a\) killing a nonzero cyclic submodule of \(Ap^*\)
kills \(V^*\).

Put \(H=Ap^*\subseteq V^*\). Every nonzero submodule of \(H\) has
annihilator \(\Ann_A(V^*)\), since it contains some nonzero \(rp^*\)
and its cyclic submodule. Choose \(G\) maximal among submodules of
\(V^*\) disjoint from \(H\), by Zorn's lemma. If \(L/G\) is a nonzero
submodule of \(V^*/G\), then \(L\) strictly contains \(G\).
By maximality \(L\cap H\ne0\), and it embeds in \(L/G\) because
\(G\cap H=0\). Its annihilator is \(\Ann_A(V^*)\), which already
kills the whole quotient. Hence \(L/G\) has that same annihilator,
so \(V^*/G\) is prime. The element \(p^*\in H\) survives.
\end{proof}

\begin{proof}[Proof of Theorem~\ref{thm:separation}]
By Lemmas~\ref{lem:ultrapower} and~\ref{lem:prime-ultrapower-quotient},
the composite \(M\xrightarrow{f} V^*\to V^*/G\) sends \(q\) to a
nonzero element. Its image is a nonzero submodule of a prime module
and hence prime. This image is the required quotient of \(M\).
\end{proof}

\begin{proof}[Proof of Corollary~\ref{cor:separation}]
Let \(S\) be a semiprime submodule of an \(A\)-module \(W\). For each
nonzero class in \(W/S\), the theorem gives a prime quotient of \(W/S\)
in which that class survives. The kernels of the composite maps from
\(W\) are prime submodules containing \(S\). Every element outside
\(S\) is omitted by one of these kernels, so their intersection is
\(S\), with the empty intersection understood as \(W\).
For left ideals, take \(W={}_AA\) and apply Lemma~\ref{lem:cyclic}.
\end{proof}

Countability is used to enumerate \(A_0\) in Lemma~\ref{lem:countable}
and to have only countably many finite systems in
Lemma~\ref{lem:tables}. The finite tests themselves use only the
finite word supports of the two operators under consideration.

\begin{remark}[Semiprime modules over commutative rings]
The condition \(r^2m=0\Rightarrow rm=0\) for all \(r,m\) is classical
semiprimeness over a commutative ring. The definition used in
Theorem~\ref{thm:separation} implies this condition, by applying it to
\(J=(r)\) and \(rm\) when \(r^2m=0\).

For the distinction, consider the example of Jenkins and Smith
\cite[Theorem~11]{JenkinsSmith}. Let \(A=k[x,y]\) for a field \(k\),
\(\mathfrak m=(x,y)\), \(v=(x,y)\in A^2\) and \(N=\mathfrak m v\).
The map \((a,b)\mapsto ya-xb\) gives \(A^2/Av\cong\mathfrak m\),
which is torsion-free, while \(Av/N\cong A/\mathfrak m\).
Thus the annihilator of every nonzero submodule of \(A^2/N\) is zero
or \(\mathfrak m\). Both ideals are prime, so \(A^2/N\) is classically
prime and hence classically semiprime. In the sense of
Theorem~\ref{thm:separation}, semiprimeness fails because
\(v+N\) is nonzero in \(\mathfrak m(A^2/N)\) and is killed by
\(\mathfrak m\). Indeed, every prime quotient of \(A^2/N\) kills
\(v+N\). If its image were nonzero, \(x\) and \(y\) would kill its
cyclic submodule and hence the whole quotient. The image of
\(v=x(1,0)+y(0,1)\) would then be zero.
\end{remark}

\subsection{The threshold for free algebras}\label{sec:threshold}

\begin{proof}[Proof of Corollary~\ref{dichotomy}]
For countable \(X\), apply Corollary~\ref{cor:separation}. Let \(X\)
be uncountable. Choose a maximal ideal \(\mathfrak m\) of \(k\), put
\(\kappa=k/\mathfrak m\), and let
\(\varphi\colon k\langle X\rangle\to\kappa\langle X\rangle\) be the
induced surjection. Take \(I\) and \(t\) from Theorem~\ref{main} for
\(\kappa\langle X\rangle\), and put \(\tilde I=\varphi^{-1}(I)\). If
\(a\,k\langle X\rangle\,a\subseteq\tilde I\), then
\(\varphi(a)\kappa\langle X\rangle\varphi(a)\subseteq I\), so
\(\varphi(a)\in I\) and \(a\in\tilde I\). Thus \(\tilde I\) is
semiprime. A prime left ideal \(\tilde P\supseteq\tilde I\) contains
\(\ker\varphi\). Hence \(\tilde P=\varphi^{-1}(\varphi(\tilde P))\),
and \(\varphi(\tilde P)\) is a prime left ideal containing \(I\). It
contains \(t\), so \(\tilde P\) contains the free generator of
\(k\langle X\rangle\) that maps to \(t\). This generator does not lie
in \(\tilde I\).
\end{proof}

Consequently, over every fixed nonzero commutative ring, the least
possible cardinality of a set of generators of an algebra that fails
prime separation is \(\aleph_1\).

\section{Other notions of prime and semiprime left ideals}
\label{sec:historical}

\subsection{Definitions}\label{sec:four-tests}
Throughout this section \(R\) denotes a ring, and \(k\) denotes a field
in the examples. For a left ideal \(I\) of \(R\), write
\(\mathbb I(I)=\{b\in R:Ib\subseteq I\}\) for its idealizer.

\begin{definition}\label{def:types}
Let \(P\) be a proper left ideal and \(I\) a left ideal of \(R\). In the
table, \(a,b\in R\), \(B,C\) are two-sided ideals of \(R\), and \(U\) is
a left ideal containing \(P\) in the middle column and containing \(I\)
in the right column. The left ideal \(P\) is \emph{extremely prime},
\emph{prime}, \emph{completely prime}, \emph{weakly prime} or
\emph{classically prime} if the implication in the corresponding row of
the middle column holds for all such \(a,b,B,C,U\). The
\emph{semiprimeness condition} in each row is the condition in the
right column, obtained by setting \(b=a\) in (E)--(W) or \(C=B\) in
(Cl) and replacing \(P\) by \(I\). When the condition holds for every
choice of the indicated variables, \(I\) is called extremely
semiprime, semiprime, weakly semiprime or classically semiprime in
rows (E), (P), (W) and (Cl), respectively. In row (C), the condition
says that \(R/I\) is endo-reduced.
{\small
\[
\begin{array}{@{}l@{\enspace}l@{\quad}l@{\quad}l@{}}
 & & \text{condition on } P & \text{condition on } I\\
\hline
\text{(E)} & \text{extremely (semi)prime}
 & ab\in P\Rightarrow a\in P\text{ or }b\in P
 & a^2\in I\Rightarrow a\in I\\
\text{(P)} & \text{(semi)prime}
 & aRb\subseteq P\Rightarrow a\in P\text{ or }b\in P
 & aRa\subseteq I\Rightarrow a\in I\\
\text{(C)} & \text{completely prime}
 & ab\in P,\ b\in\mathbb I(P)\Rightarrow a\in P\text{ or }b\in P
 & a^2\in I,\ a\in\mathbb I(I)\Rightarrow a\in I\\
\text{(W)} & \text{weakly (semi)prime}
 & aRb\subseteq P,\ Rb\subseteq\mathbb I(P)\Rightarrow a\in P\text{ or }b\in P
 & aRa\subseteq I,\ Ra\subseteq\mathbb I(I)\Rightarrow a\in I\\
\text{(Cl)} & \text{classically (semi)prime}
 & BCU\subseteq P\Rightarrow BU\subseteq P\text{ or }CU\subseteq P
 & B^2U\subseteq I\Rightarrow BU\subseteq I
\end{array}
\]}%
\end{definition}

The term \emph{extremely prime} is used in \cite[Section~2.2]{AryapoorPrime},
which also notes the use of \emph{completely prime} for (E).
T\"orner and Zima call the right-sided form of the semiprimeness
condition in (E) \emph{complete semiprimeness} \cite[Section~1]{ToernerZima}. Koh called
the right-sided form of (P) \emph{prime type} \cite{KohPrimeType}.
Condition (C) is the left-sided form of Reyes's completely prime right
ideals \cite[Definition~2.1]{Reyes}. Right multiplication identifies
\(\operatorname{End}_R(R/I)\) with the opposite ring of \(\mathbb I(I)/I\).
Hence \(P\) satisfies (C) if and only if every nonzero endomorphism of
\(R/P\) is injective \cite[Proposition~2.5]{Reyes}. The condition on \(I\)
in row (C) says that \(R/I\) is \emph{endo-reduced}, meaning that
\(\operatorname{End}_R(R/I)\) is reduced \cite{KimuliSsevviiri}.

The idealwise forms of (P) and (W) differ only in the left ideals that
are tested. A proper left ideal \(P\) is prime if and only if
\(UV\subseteq P\) implies \(U\subseteq P\) or \(V\subseteq P\) for all
left ideals \(U,V\), and \(I\) is semiprime if and only if
\(U^2\subseteq I\) implies \(U\subseteq I\) for all left ideals \(U\).
Van der Walt showed that \(P\) is weakly prime if and only if the same
holds for all left ideals \(U,V\supseteq P\), that is, \(UV\subseteq P\)
implies \(U=P\) or \(V=P\) \cite[Proposition~1.1]{VanDerWalt}. Since
\((I+Ra)(I+Rb)\subseteq I\) holds exactly when \(aRb\subseteq I\) and
\(IRb\subseteq I\), in the same way \(I\) is weakly semiprime if and only
if \(U\supseteq I\) and \(U^2\subseteq I\) imply \(U=I\).
The two conditions in (Cl) can also be tested on elements. One may
take \(U=P+Ra\), respectively \(U=I+Ra\): \(P\) is classically prime if
and only if \(BCa\subseteq P\) implies \(Ba\subseteq P\) or
\(Ca\subseteq P\), and \(I\) is classically semiprime if and only if
\(B^2a\subseteq I\) implies \(Ba\subseteq I\). (If \(BCU\subseteq P\),
\(Ba\not\subseteq P\) and \(Ca'\not\subseteq P\) with \(a,a'\in U\), the
elementwise condition for \(a\), \(a'\) and \(a+a'\) gives a
contradiction.) Equivalently, \(P\) is classically prime if and only if every nonzero
submodule of \(R/P\) has a prime annihilator
\cite{Behboodi,BehboodiShojaee,GroenewaldSsevviiriClassical}.
For a two-sided ideal \(I\) we have \(\mathbb I(I)=R\), so (E) and (C)
agree with complete primeness of two-sided ideals, and (W) agrees with (P). The same identifications hold for the
corresponding semiprimeness conditions.

\begin{figure}[htbp]
\centering
\[
\begin{array}{c@{\enspace}c@{\enspace}c@{\enspace}c@{\enspace}c}
 \text{Extremely prime (E)}&\Longrightarrow&\text{Prime (P)}
 &\Longrightarrow&\text{Classically prime (Cl)}\\
 \Downarrow&&\Downarrow&&\\
 \text{Completely prime (C)}&\Longrightarrow&\text{Weakly prime (W)}&&
\end{array}
\]
\caption{Implications among five notions of a prime left ideal.}
\label{fig:prime-implications}
\end{figure}

Each condition on \(P\) implies its semiprimeness condition for \(I=P\).
For (P)\(\Rightarrow\)(Cl), Lemma~\ref{lem:cyclic} makes \(R/P\) a prime
module. Every nonzero submodule therefore has annihilator
\(\Ann_R(R/P)\), which is the prime two-sided core of \(P\) by
Lemma~\ref{core}. The other arrows follow directly from the tests
in Definition~\ref{def:types}.

The arrows in Figure~\ref{fig:prime-implications} and their composites
give all general implications among the prime notions. The zero ideal of \(M_2(k)\) satisfies (P)
and fails (C). In Example~\ref{ex:triangular}, \(P\) satisfies (C) and
fails (P) and (Cl). For a weakly prime left ideal that fails both
(P) and (C), specialize van der Walt's construction
\cite[Example~1.2]{VanDerWalt} to \(S=k[\varepsilon]/(\varepsilon^2)\).
Put \(R=M_2(S)\) and \(L=Re_{11}\). Since \(LR=R\), the
condition \(LRb\subseteq L\) forces \(b\in L\), so \(L\) is weakly
prime. The element \(e=\varepsilon 1_R\) is central, lies outside \(L\),
and satisfies \(e^2=0\) and \(eRe=0\). The left ideal in
Section~\ref{sec:classical} satisfies (Cl) and fails (W).

\subsection{Intersections}\label{sec:intersections}
Table~\ref{tab:tests} compares the semiprimeness condition in each row
of Definition~\ref{def:types} with being an intersection of left ideals
that satisfy the corresponding primeness condition. For (E), (P) and
(Cl), the hypothesis of the semiprimeness condition passes from the
intersection to each member, with \(U+P\) in place of \(U\) for (Cl). For (C) and
(W), the idealizer condition can fail to pass to the members
(Example~\ref{ex:triangular}).

For two-sided \(I\), the semiprimeness conditions in (E) and (C) both
express reducedness of \(R/I\). By Andrunakievich and Ryabukhin's theorem,
such an \(I\) is the intersection of the completely prime two-sided
ideals containing it \cite[Theorem~1]{AndrunakievichRyabukhin}. These
two-sided factors satisfy both (E) and (C).

\begin{table}[htbp]
\centering
\begin{tabular}{@{}lll@{}}
\hline
Type & Intersections satisfy & Semiprimeness condition implies\\
 & the semiprimeness condition & being an intersection\\
\hline
(E) & yes & no (Example~\ref{ex:extreme})\\
(P) & yes & no (Theorem~\ref{main}); yes for countably\\
 & & generated algebras (Corollary~\ref{cor:separation})\\
(C) & no (Example~\ref{ex:triangular}) & no (Example~\ref{ex:reyes});
 yes if \(R/I\) has\\
 & & finite length (Proposition~\ref{prop:reyes-length})\\
(W) & no (Example~\ref{ex:triangular}) & yes (Proposition~\ref{prop:weak-separation})\\
(Cl) & yes & no (Example~\ref{ex:classical}); yes in the\\
 & & cases of Section~\ref{sec:classical}\\
\hline
\end{tabular}\par\smallskip
\caption{Intersection properties of the five types in Definition~\ref{def:types}.}
\label{tab:tests}
\end{table}

\begin{example}\label{ex:triangular}
In \(R=T_2(k)\), put \(P=ke_{11}\) and \(Q=ke_{12}+ke_{22}\).
The left ideal \(Q\) is maximal, and \(R/P\cong Q\) has endomorphism
ring \(k\), so both are completely prime and hence weakly prime. Their
intersection \(0\) fails the semiprimeness conditions in (C) and (W),
since \(\mathbb I(0)=R\) and \(a=e_{12}\) satisfies \(a^2=0\) and
\(aRa=0\). As \(a\notin P\), the left ideal \(P\) is not semiprime and
fails (P). The module \(Q\) is
faithful, while \(R\) has the nonzero square-zero ideal \(ke_{12}\), so
\(P\) fails (Cl).
\end{example}

\subsection{Weakly prime left ideals}\label{sec:weak-primeness}
Every semiprime left ideal is weakly semiprime, so the following
proposition contains the consequence of Van der Walt's work quoted in
the introduction.

\begin{proposition}\label{prop:weak-separation}
Every weakly semiprime left ideal is the intersection of the weakly prime
left ideals containing it.
\end{proposition}

\begin{proof}
We use the idealwise forms of (W). Let \(I\) be weakly semiprime and
\(a_0\notin I\). Following Van der Walt's
sequence argument \cite[proof of Proposition~3.2]{VanDerWalt}, choose
\(a_{n+1}\in (I+Ra_n)^2\setminus I\) at each step. By Zorn's lemma,
choose \(P\) maximal among the left ideals containing \(I\) and avoiding
the sequence. A left ideal properly containing \(P\) contains some \(a_n\)
and the entire tail, since \(a_{n+1}\in I+Ra_n\). Thus any two left ideals
\(U,V\) properly containing \(P\) contain a common \(a_n\) and satisfy
\(a_{n+1}\in (I+Ra_n)^2\subseteq UV\). Hence \(UV\not\subseteq P\),
so \(P\) is weakly prime and omits \(a_0\).
\end{proof}

By Example~\ref{ex:triangular}, an intersection of weakly prime left
ideals need not be weakly semiprime. In van der Walt's criterion, a
left ideal maximal among those disjoint from an \(m\)-system is prime
if its two-sided core is maximal among the two-sided ideals disjoint
from the same system \cite[Proposition~1.6]{VanDerWalt}.

\subsection{Completely prime left ideals}\label{sec:reyes}
\begin{proposition}\label{prop:reyes-length}
If \(R/I\) has finite length and is endo-reduced,
then \(I\) is a finite intersection of completely prime left ideals.
\end{proposition}

\begin{proof}
Decompose \(R/I\) as a finite direct sum of indecomposable modules \(M_j\).
By Fitting's lemma, an endomorphism of \(M_j\) is either invertible or
nilpotent. A nilpotent endomorphism extends by zero to \(R/I\), so
reducedness of \(\operatorname{End}_R(R/I)\) forces it to vanish. Hence
the kernel of each projection \(R\to M_j\) is completely prime,
and these kernels intersect in \(I\).
\end{proof}

\begin{example}\label{ex:reyes}
Let \(D=k[t]\), \(R=T_2(D)\) and \(I=De_{11}+t^2De_{12}\). A matrix lies
in \(\mathbb I(I)\) exactly when its upper right entry is divisible by
\(t^2\), and the lower right entry gives \(\mathbb I(I)/I\cong D\).
Thus \(R/I\) is endo-reduced. Let \(P\supseteq I\) be
completely prime. Since \(t1_R\) is central and
\((te_{12})(t1_R)\in I\), we get \(te_{12}\in P\) or \(t1_R\in P\), and
in the latter case \(te_{12}=e_{12}(t1_R)\in P\). So every such \(P\)
contains \(te_{12}\notin I\).
\end{example}

Reyes's result for Ore sets implies that a central element \(z\) lies
in every completely prime left ideal containing \(I\) exactly when
\(z^n\in I\) for some \(n\ge1\) \cite[Example~5.19]{Reyes}.
For polynomial rings in central variables over the quaternions, Alon
and Paran identified this intersection with the ideal of polynomials
vanishing at the common zeros with commuting coordinates
\cite[Theorem~1.2]{AlonParan}.

\subsection{Extremely prime left ideals}\label{sec:extreme}
\begin{example}\label{ex:extreme}
Let \(R=k\langle r,a,v\rangle/(va-r^2)\) and \(I=Rra\). We show that
\(I\) is extremely semiprime, but \(r^2\notin I\) lies in every
extremely prime \(P\supseteq I\): such \(P\) contains \(r\) or \(a\),
hence \(r^2=va\).

Occurrences of \(va\) are disjoint, and replacing one by \(rr\)
creates no new occurrence. Each word therefore has a unique normal
form. The linear normal-form map kills every \(u(va-r^2)w\), for words
\(u,w\), so the words avoiding \(va\) form a basis. The left ideal \(I\)
is spanned by those ending in \(ra\). In particular \(r^2\notin I\). Order these words by length
and then lexicographically with \(v<r<a\). The only rewritings are
\(v\cdot a=rr\) at a junction, and \(rr\) lies after \(vv,vr\) and before
\(ra,aa\). Hence multiplication by a letter on either side preserves
the order, and the largest word of a product of two nonzero elements is the
normal form of the product of their largest words.

Suppose \(f^2\in I\) and \(f\notin I\), and let \(f_0\ne0\) be the part
of \(f\) supported outside \(I\). Then \(ff_0\in I\). Let \(g,h\) be the
largest words of \(f,f_0\). The normal form of \(gh\) ends in \(ra\),
while \(h\) does not. Thus \(h=1\), or \(h=a\) and \(g\) ends in \(r\),
or \(h=a^2\) and \(g\) ends in \(v\). If \(h=1\), then \(f_0=c\)
is a nonzero scalar, so \(f=c^{-1}ff_0\in I\). Otherwise \(g\) does not end in \(ra\),
so \(g\) occurs in \(f_0\) and \(g=h\), contrary to the ending of \(g\).
\end{example}

\subsection{Classically prime left ideals}\label{sec:classical}
For commutative rings, classical primeness and classical semiprimeness
of ideals agree with the ordinary notions. Behboodi proved classical
prime separation for classically semiprime submodules over fully
bounded Noetherian prime rings of classical Krull dimension at most
one \cite{Behboodi}. Taking the regular left module gives the
corresponding result for left ideals.

To see that (Cl) does not imply (W), let \(D=k[t]\),
\(R=M_2(D)\), and \(I=tRe_{22}\). As a left module,
\(R/I\cong D^2\oplus(D/tD)^2\). A submodule with nonzero projection
to \(D^2\) has annihilator zero, since its projection contains \(fD^2\)
for some \(0\ne f\in D\) by the matrix units. Any other nonzero submodule
has annihilator \(tR\). Both annihilators are prime, so \(I\) is
classically prime. But \(a=t1_R\) and \(b=e_{22}\) are outside \(I\) and
satisfy \(aRb\subseteq I\) and \(IRb\subseteq I\).

\begin{example}\label{ex:classical}
Let \(D\) be the first Weyl algebra over \(\mathbb Q\), with nonunit
generator \(x\). Put \(S=xD\), \(R=\mathbb Q(1,1)+(S\times S)\),
and \(I=R(x,-x)\). Since \(SdS=S\) for \(d\ne0\) in the simple domain
\(D\), an ideal of \(R\) with nonzero first coordinate contains
\(A=S\times0\). The second coordinate gives \(B=0\times S\) in the
same way. Since \(R/(A+B)\cong\mathbb Q\), all two-sided ideals are
\(0,A,B,A+B,R\). They are idempotent, so \(I\) is classically semiprime.

If \(P\supseteq I\) is classically prime, \(AB=0\) forces \(A\subseteq P\)
or \(B\subseteq P\), by primeness of \(\Ann_R(R/P)\). In either case
\((x,0)\in P\), since \((x,0)-(0,x)\in I\). But \((x,0)\notin I\):
the equation \((a,b)(x,-x)=(x,0)\) would force \((a,b)=(1,0)\notin R\).
\end{example}

\section*{Acknowledgment}
The author used Claude (Anthropic) and ChatGPT (OpenAI) extensively
to explore proof strategies and revise the exposition. All mathematical
results were human-verified. The author takes responsibility for the
content of the manuscript.

\end{document}